\documentclass{amsart}
\usepackage{amssymb,euscript,amsmath, mathrsfs}

\newcounter{ENUM}

\newcommand{\margh}[1]{}

\input xy
\xyoption{all}
\CompileMatrices

\def\sE{{\mathscr E}}
\def\sF{{\mathscr F}}

\def\sK{{\mathscr K}}
\def\sL{{\mathscr L}}

\def\sO{{\mathscr O}}
\def\sQ{{\mathscr Q}}

\def\Supp{\operatorname{Supp}}
\newcommand{\el}{$\ell$}

\newtheorem{thm}{Theorem}[section]
\newtheorem{prop}[thm]{Proposition}

\newtheorem{cor}[thm]{Corollary}

\theoremstyle{definition}
\newtheorem{defn}[thm]{Definition}

\newtheorem{ex}[thm]{Example}

\theoremstyle{remark}
\newtheorem{rem}[thm]{Remark}

\numberwithin{equation}{section}
\numberwithin{figure}{section}

\begin{document}
\title{Stability of vector bundles on curves and degenerations}
\author{Brian Osserman}
\begin{abstract} 
We observe that if we are interested primarily in degeneration arguments,
there is a weaker notion of (semi)stability for vector bundles on 
reducible curves, which is sufficient for many applications, and does not 
depend on a choice of polarization. We introduce and explore the basic
properties of this alternate notion of (semi)stability. In a complementary
direction, we record a proof of the existence of semistable extensions
of vector bundles in suitable degenerations.
\end{abstract}

\maketitle

\section{Introduction}

Typically, when considering (semi)stability of vector bundles on a 
reducible curve $X$, one works with respect to a polarization, which
assigns weights to the different irreducible components of $X$. The
reason for this is that subsheaves of vector bundles may have different
ranks on different components of $X$, and one needs to determine how
to weigh these ranks. Put differently, the Hilbert polynomials of 
subsheaves will now depend nontrivially on the choice of an ample line
bundle. Such an approach is very natural for working with coarse moduli
spaces, but it introduces undesirable technical complications. The 
present note is based on the observation that if one is 
interested in (semi)stability in the context of degeneration techniques 
(typically applied to higher-rank Brill-Noether theory), a 
weaker definition suffices, in which one only considers constant-rank
subsheaves. The resulting definition is independent of polarization,
and both better-behaved and easier to verify than the standard one.
Although it is not well-suited to coarse moduli space constructions,
it is still an open condition in families, and hence works well in the
context of moduli stacks.
We hope that the new perspective will not only streamline proofs of
known results (see for instance Remark \ref{rem:canon-exist}), but will 
also open the door to improved existence results in higher-rank 
Brill-Noether theory. There are even hints that it may have a role to
play in specialization arguments, where \textit{a priori} one would
expect to want to make use of the usual stronger notion of semistability;
see Remark \ref{rem:petri-map}.

In a complementary direction, we record in Proposition \ref{prop:specialize}
a natural statement which does not seem to have appeared in the literature
regarding existence of semistable extensions of vector bundles 
with respect to a chosen polarization.

We assume throughout that all curves are proper, geometrically reduced and 
geometrically connected, but not necessarily irreducible.

We now introduce the notion of \el-semistability (short for ``limit
semistability'').

\begin{defn} Let $X$ be a nodal (possibly reducible) curve, and
$\sE$ a vector bundle of rank $r$ on $X$. We say that $\sE$ is
$\boldsymbol{\ell}$-\textbf{semistable} (respectively, 
$\boldsymbol{\ell}$-\textbf{stable}) if
for all proper subsheaves $\sF \subseteq \sE$ having constant rank
$r'$ on every component of $X$, we have
\begin{equation}\label{eq:defn}
\frac{\chi(\sF)}{r'} \leq \frac{\chi(\sE)}{r} \left(\text{respectively, }
\frac{\chi(\sF)}{r'} < \frac{\chi(\sE)}{r}\right).
\end{equation}
\end{defn}

Thus, if $X$ is irreducible, we recover the usual definition of
(semi)stable, while in the reducible case we have a weaker definition,
which does not involve a polarization. Note that an \el-stable
vector bundle may not be stable with respect to any polarization;
see Example \ref{ex:ell-stable-not-stable}.

Our main results are the following observations:

\begin{prop}\label{prop:open}
Both \el-semistability and \el-stability are open in families. 
\end{prop}

\begin{prop}\label{prop:twist}
Both \el-semistability and \el-stability are closed under
tensor product with line bundles.
\end{prop}

Consequently, we find:

\begin{cor}\label{cor:open}
Let $\pi:X \to S$ be a family of curves, with $S$ 
the spectrum of a DVR, such that $\pi$ has smooth generic fiber $X_{\eta}$, 
and nodal special fiber $X_0$. Let $\sE$ be a vector bundle on $X$.
If $\sE$ is \el-semistable (respectively, \el-stable) on $X_0$, then
$\sE$ is semistable (respectively, stable) on $X_{\eta}$.
\end{cor}

\begin{cor}\label{cor:twist}
Let $X$ be a nodal curve, $\sE$ a vector bundle on $X$, and
$\sL$ a line bundle on $X$. If $\sE \otimes \sL$ is semistable (respectively,
stable) with respect to some polarization on $X$, then $\sE$ is 
\el-semistable (respectively, \el-stable).
\end{cor}

Finally, we demonstrate that \el-semistability behaves very well with 
respect to gluing of subcurves.

\begin{prop}\label{prop:naive}
Let $X=Y \cup Z$ be a nodal curve, with the subcurves $Y$ and $Z$ meeting
at a single node $P$. Given a vector bundle $\sE$ on $X$ of rank $r$, if 
$\sE|_Y$ and $\sE|_Z$ are \el-semistable on $Y$ and $Z$ respectively, then 
$\sE$ is \el-semistable on $X$.
If further there do not exist subsheaves
$\sF_Y\subseteq \sE|_Y$ and $\sF_Z\subseteq \sE|_Z$ of some constant rank 
$r'$ which glue to one another at $P$, and which satisfy
$$\chi(\sF_Y)/r' = \chi(\sE|_Y)/r \quad \text{ and } \quad
\chi(\sF_Z)/r' = \chi(\sE|_Z)/r,$$
then $\sE$ is \el-stable.
\end{prop}

Here, when we say that $\sF_Y$ and $\sF_Z$ glue to one another at $P$,
we mean that in a neighborhood of $P$ they are obtained as the restrictions
to $Y$ and $Z$ respectively of a subbundle of $\sE$.

\begin{cor}\label{cor:naive}
Let $X$ be a curve of compact type, and $\sE$ a vector bundle
on $X$. If $\sE|_Y$ is semistable for every component $Y$ of $X$, then
$\sE$ is \el-semistable. If moreover
there does not exist a vector subbundle $\sF \subseteq \sE$ which is
weakly destabilizing on every component of $X$,
then $\sE$ is \el-stable.
\end{cor}

Here, a subbundle $\sF \subseteq \sE|_Y$ is \textbf{weakly destabilizing}
if $\chi(\sF)/r'=\chi(\sE|_Y)/r$, where $r'$ and $r$ are the ranks of $\sF$
and $\sE$, respectively. Corollary \ref{cor:naive} may be deduced from
Corollary \ref{cor:twist} and analogous results in the literature on
usual stability (see for instance Proposition 1.2 of \cite{te7}),
but there are many \el-semistable vector bundles
which are not of the form considered in Corollary \ref{cor:naive},
and Proposition \ref{prop:naive} provides a powerful tool for building
them one block at a time. We
hope that systematically considering such bundles will lead to better
existence results in higher-rank Brill-Noether theory (for instance, in
the direction of the Bertram-Feinberg-Mukai conjecture) than those which 
have been obtained to date.

\section{Proofs}

We now give the proofs of the claimed results.

\begin{proof}[Proof of Proposition \ref{prop:open}]
The main observation is that if $\pi:X \to S$ is a (flat, proper) family 
of curves
where $S$ is connected and locally Noetherian, and $\sE$ is a coherent sheaf 
on $X$, flat over $S$, if there exists $s \in S$ such that $\sE$ has the same
rank $r$ generically on every component of the fiber $X_s$, then the same
is true for all $s \in S$. It clearly suffices to handle the case that $S$
is irreducible, so in this case, we first prove that the statement holds
for the generic point $\eta$ of $S$, and then for all $s' \in S$. Since
the hypotheses are preserved under base change and the conclusion may be
tested after base change, we therefore reduce to the case that $S$ is
the spectrum of a DVR, and we wish to prove that $\sE$ has rank $r$ on
every component of the generic fiber $X_{\eta}$ if and only it has rank
$r$ on every component of the special fiber $X_0$. 
However, by flatness over $S$, the open subset of $X$ on which $\sE$ is 
locally free must meet every component of $X_0$: indeed, the support of 
any torsion of $\sE$ cannot contain any generic point of the special fiber 
without creating torsion over $S$.
Since every component of
$X_{\eta}$ must contain at least one component of $X_0$ in its closure,
and every component of $X_0$ is in the closure of some component of 
$X_{\eta}$,
the desired statement follows.

It thus follows that the locus in a given Quot scheme which consists 
of quotient sheaves having equal rank on each component is a union of
connected components of the Quot scheme, and is in particular proper.
The proposition then follows from the usual argument for openness of
(semi)stability (see for instance Proposition 2.3.1 of \cite{h-l}).
\end{proof}

\begin{proof}[Proof of Proposition \ref{prop:twist}]
Let $\sE$ be a vector bundle of rank $r$ on a curve $X$, and $\sL$ a line 
bundle. Observe that, for any $r'<r$, tensoring with $\sL$ induces a 
bijection between subsheaves of $\sE$ of pure rank $r'$ and subsheaves
of $\sE \otimes \sL$ of pure rank $r'$. It is thus enough to observe that
for any subsheaf $\sF \subseteq \sE$ of pure rank $r'$, we have
$$\frac{\chi(\sE)}{r}-\frac{\chi(\sF)}{r'}
=\frac{\chi(\sE \otimes \sL)}{r}-\frac{\chi(\sF\otimes \sL)}{r'},$$
which we prove by showing that
$\chi(\sF \otimes \sL)=\chi(\sF)+r'\deg \sL$, and
$\chi(\sE \otimes \sL)=\chi(\sE)+r \deg \sL$. This is presumably standard,
but since the proof contains minor subtleties in the case of a reducible 
curve, we include it for the sake of completeness.

Let $D$ be a sufficiently ample effective divisor supported on
the smooth locus of $X$ such that $\sL(D)$ has a section $s$ which is 
nonvanishing at the nodes of $X$.
Then using the short exact sequences induced by $s$ gives us that
$\chi(\sF \otimes \sL(D))=\chi(\sF)+r'(\deg \sL + \deg D)$
and
$\chi(\sE \otimes \sL(D))=\chi(\sE)+r(\deg \sL + \deg D)$.
Then, the exact sequences induced by the canonical inclusion
$\sL \hookrightarrow \sL(D)$ yields the desired identities.
\end{proof}

Note that Corollary \ref{cor:open} is an immediate consequence of
Proposition \ref{prop:open}, and Corollary \ref{cor:twist} is an
immediate consequence of Proposition \ref{prop:twist}.

\begin{proof}[Proof of Proposition \ref{prop:naive}]
Let $\sF$ be a subsheaf of $\sE$ of
constant rank $r'$. Let $\sF_Y$ (respectively, $\sF_Z$) denote
the quotient of $\sF|_Y$ (respectively, $\sF|_Z$) by its torsion subsheaf, 
or equivalently, the image of $\sF|_Y$ inside $\sE|_Y$ (respectively, of
$\sF|_Z$ inside $\sE|_Z$). Then there is an integer $r_P$ between $0$ 
and $r'$ described as follows: if $\sQ$ is the cokernel of
$$\sF \hookrightarrow \sF_Y \oplus \sF_Z,$$
then one checks that injectivity is preserved after restriction to $P$,
so we have an induced exact sequence 
$$0 \to \sF|_P \to \sF_Y|_P \oplus \sF_Z|_P \to \sQ|_P \to 0,$$
and we let $r_P$ be the dimension of $\sQ|_P$. 
We then have that 
$r_P=r'$ if and only if $\sF$ is locally free at $P$,
and furthermore, if $\widetilde{\sF}_Y$ denotes the saturation of $\sF_Y$ 
at $P$ in $\sE|_Y$, and similarly for $\widetilde{\sF}_Z$, then $r_P$ has 
the property
that the quotients $\widetilde{\sF}_Y/\sF_Y$ and $\widetilde{\sF}_Z/\sF_Z$
each have dimension at least $r'-r_P$ at $P$.
Now, we carry out the following calculation:
\begin{align*}
\frac{\chi(\sF)}{r'} & \leq \frac{1}{r'}\left(\chi(\sF_Y) + \chi (\sF_Z)
- r_P\right) \\
& \leq \frac{1}{r'}\left(\chi(\widetilde{\sF}_Y)+\chi(\widetilde{\sF}_Z)
- 2(r'- r_P) - r_P\right) \\
& \leq \frac{1}{r}\left(\chi(\sE|_Y)+\chi(\sE|_Z)\right)
- \frac{2r'-r_P}{r'} \\
& = \frac{\chi(\sE)}{r}+ 1- \frac{2r'-r_P}{r'} \\
& = \frac{\chi(\sE)}{r}- \frac{r'-r_P}{r'} \\
& \leq \frac{\chi(\sE)}{r}.
\end{align*}
Thus, we get that $\sE$ is \el-semistable, and is in fact \el-stable unless
there exists some $\sF$ with
$\chi(\widetilde{\sF}_Y)/r'=\chi(\sE|_Y)/r$,
$\chi(\widetilde{\sF}_Z)/r'=\chi(\sE|_Z)/r$,
and $r_P=r'$. The condition $r_P=r'$ implies that 
$\sF_Y=\widetilde{\sF}_Y$ and $\sF_Z=\widetilde{\sF}_Z$, and that 
$\sF_Y$ must glue to $\sF_{Z}$ at $P$, as desired.
\end{proof}

Corollary \ref{cor:naive} follows immediately from Proposition 
\ref{prop:naive} by induction on the number of components of $X$.

\section{Further discussion}

It is instructive to compare \el-semistability to usual semistability
in the case of rank 2 and degree $2g-2$. This case is the
subject of the Bertram-Feinberg-Mukai conjecture, and has consequently
received a great deal of attention. A vector bundle $\sE$ of rank $2$
and degree $2g-2$ has $\chi(\sE)=0$, so we see that although
the usual definition of (semi)stability calls for a polarization, the
resulting definition in fact does not depend on the polarization. This
therefore presents a natural context in which to compare the definitions.

\begin{ex}\label{ex:ell-stable-not-stable} Let $X$ be a chain of smooth
curves $Y_1,\dots,Y_n$, glued together at nodes. Let $\sE$ be a vector 
bundle of rank $2$ and degree $2g-2$, satisfying the condition for
\el-stability of Corollary \ref{cor:naive}. Suppose further that
$d_1:=\deg \sE|_{Y_1}$ is even, and $\sE|_{Y_1}$ is strictly semistable. Then
we see that even though $\sE$ is \el-stable, it is not stable on $X$.
Indeed, let $g_1$ be the genus of $Y_1$; then if $d_1 \geq 2g_1$, the 
condition for stability is violated by the subsheaf of $\sE$ consisting of 
sections which vanish on the complement of $Y_1$, 
while if $d_1 \leq 2g_1-2$, the condition for stability is violated by 
the subsheaf of $\sE$ consisting of sections which vanish on $Y_1$.
\end{ex}

The peculiarities of the case $\chi(\sE)=0$ lead to subtleties in certain 
aspects of Teixidor i Bigas' \cite{te5} and \cite{te2}. These subtleties
are addressed by twisting arguments, and we take the opportunity to 
discuss how they fit into the context of \el-stability.

\begin{rem}\label{rem:petri-map} In \cite{te2}, a key point is to place
suitable conditions on the vector bundle $\sE_0$ obtained as a specialization
from a semistable vector bundle on the smooth generic fiber; this is
carried out in Claim 2.3. However, Claim 2.3 does not apply directly to
the case of interest, because when $\chi=0$ it is not possible to choose a 
polarization with the required non-integrality property. Instead, as
described at the beginning of \S 3 of \cite{te2}, one uses twisting to
complete the argument, as follows.
Let $\pi:X \to S$ be the family of curves used for
the degeneration and let $\sE_{\eta}$ be a vector bundle on the generic
fiber $X_{\eta}$, with canonical determinant. Choose $D$ any divisor on 
$X$ of non-zero relative degree; rather than extending $\sE_{\eta}$ right 
away, we instead twist by $D$, and then extend $\sE_{\eta}(D|_{X_{\eta}})$ 
to a bundle $\sE'$ semistable with respect to a polarization satisfying the 
stated condition.\footnote{In fact, the argument in \cite{te2} is slightly
more complicated, but can be simplified to the above using Proposition
\ref{prop:specialize}.}
Then Claim 2.3 of \textit{loc.\ cit.}\ shows that $\sE'$ has the desired
properties on each irreducible component and in a neighborhood of each
node, and it follows that if we set $\sE_0=\sE'(-D)|_{X_0}$, we obtain an 
extension of the original $\sE_{\eta}$ which has the desired properties.

For us, the relevant point is that, because semistability is not preserved 
by twisting, there is no reason to think that $\sE_0$ is semistable, but at 
least it follows from Corollary \ref{cor:twist} that it is \el-semistable. 
This hints that even though for specialization arguments it is natural to
try to take advantage of the stronger properties afforded by semistability
(with respect to a polarization), there may also be a role for 
\el-semistable vector bundles.
\end{rem}

\begin{rem}\label{rem:canon-exist} In \cite{te5}, in order to stay within
the framework of usual stability, one needs to make an
argument similar to that of Remark \ref{rem:petri-map}, because the 
situation is precisely
as in Example \ref{ex:ell-stable-not-stable}, so the natural underlying
vector bundles are in fact not stable under any choice of polarization.
But if $\pi:X \to S$ and $D$ are as in the previous remark,
and $\sE$ is a vector bundle on $X$ underlying the relevant limit
linear series and its smoothing, we show that $\sE$ is stable
on the generic fiber $X_{\eta}$ as follows. The twisted bundle $\sE(D)$ has 
nonzero Euler characteristic, so after a possible further twist to 
redistribute degrees, it is stable on $X_0$ for a suitable choice
of polarization. It thus follows that $\sE(D)|_{X_{\eta}}$ is stable, and
hence, since $X_{\eta}$ is smooth, that $\sE|_{X(\eta)}$ is also stable,
as desired. 

Although the above argument works, it seems much simpler to argue that for 
the original $\sE$, although $\sE|_{X_0}$ is strictly semistable, it is 
\el-stable by Proposition \ref{prop:naive}, and therefore $\sE|_{X_{\eta}}$ 
is stable, as needed.
\end{rem}

\section{Semistable extensions}

Although the following result on specialization of vector bundles
under degeneration is a straightforward application of standard
techniques, it does not appear to be stated anywhere in the literature.
Because it complements the main subject of the present note, we take
the opportunity to record its proof.

\begin{prop}\label{prop:specialize} Let $\pi:X \to B$ be a flat,
proper morphism with $B$ the spectrum of a DVR, generic fiber $X_{\eta}$
a smooth curve, and special fiber $X_0$ a nodal curve. Suppose that
$X$ is regular. Then for any polarization $w$ on $X_0$, and any
semistable vector bundle $\sE_{\eta}$ on $X_{\eta}$, there exists a
vector bundle $\sE$ on $X$ such that $\sE|_{X_{\eta}} \cong \sE_{\eta}$
and $\sE|_{X_0}$ is semistable with respect to $w$.
\end{prop}

Recall that a polarization $w$ is a positive rational weighting of the
components of $X_0$ adding to $1$. Semistability with respect to such a
polarization is equivalent to semistability with respect to an ample
divisor supported on the smooth locus of $X_0$ -- clearing denominators in 
$w$ describes the distribution of degrees of the divisor in question.

\begin{proof} 
First, recall that a reflexive sheaf on a regular $2$-dimensional scheme
is necessarily locally free.
Thus, by extending $\sE_{\eta}$
to any coherent sheaf on $X$, and then taking the reflexive hull, we 
obtain a vector bundle $\sE'$ on $X$ extending $\sE_{\eta}$. It remains
to show that the desired $\sE$ can be realized as a subsheaf of $\sE'$.
This follows the standard argument of Langton (see Theorem 2.1.B of
\cite{h-l}); all that needs to be checked is that the subsheaves 
considered inductively in the argument in question remain locally free
at each step. But these subsheaves are obtained by considering the kernel 
$\sK$ of composed maps of the form
$$\sE \twoheadrightarrow \sE|_{X_0} \twoheadrightarrow \sF,$$
where $\sE$ is a vector bundle and
$\sF$ is the quotient sheaf corresponding to a maximally
destabilizing subbundle of $\sE|_{X_0}$. In particular, the kernel of
$\sE|_{X_0} \twoheadrightarrow \sF$ is saturated, so $\sF$ is pure of 
dimension $1$. Now, purity implies that for any closed point 
$x \in \Supp(\sF)$,
the stalk $\sF_x$ has depth $1$, so it follows 
from the Auslander-Buchsbaum formula that the projective dimension of
$\sF_x$ is also $1$. Using the Tor exact sequence, we conclude that $\sK_x$
is flat over $\sO_{X,x}$, and hence $\sK$ is locally free, as desired.
\end{proof}

If we drop the regularity hypothesis on $X$, then it is always possible to
blow up $X$ at nodes of $X_0$ to resolve any singularities; this only
introduces chains of rational components at the nodes of $X_0$, and then
it follows from Proposition \ref{prop:specialize} that we can extend any
vector bundle while preserving semistability.

\bibliographystyle{amsalpha}
\bibliography{gen}

\end{document}